\documentclass[preprint,12pt]{elsarticle}

\usepackage{amsfonts,latexsym,amsmath,graphicx,psfrag}
\usepackage{amssymb}

\newtheorem{theorem}{Theorem}

\newtheorem{lemma}[theorem]{Lemma}
\newtheorem{proposition}[theorem]{Proposition}
\newtheorem{definition}[theorem]{Definition}
\newproof{proof}{Proof}

\newcommand{\notmid}{\mid\kern-0.5em\not\kern0.5em}

\newcommand{\mnewtext}{}  
\newcommand{\mltext}{}      
\newcommand{\mtext}{}

\newcommand{\HOX}[1]{\marginpar{\footnotesize #1}}

\newcommand{\dd}{\mathrm{d}}

\newcommand{\R}{\mathbb{R}}
\newcommand{\C}{\mathbb{C}}

\def\Z{{\mathbb Z}}

\newcommand{\ignore}[1]{}

\newcommand{\mstrut}[1]{\mbox{\rule{0mm}{#1}}}

\newcommand{\barr}{\begin{array}}
\newcommand{\earr}{\end{array}}

\def\bfo{\begin{eqnarray*}}
\def\efo{\end{eqnarray*}}
\def\ba{\begin{eqnarray*}}
\def\ea{\end{eqnarray*}}
\def\beq{\begin{eqnarray}}
\def\eeq{\end{eqnarray}}

\def\hat{\widehat}
\def\tilde{\widetilde}

\def\T{{\mathcal T}}

\newcommand{\V}{{\cal V}} 
\newcommand{\W}{{\cal W}}

\def\bra{\langle}
\def\cet{\rangle}
\def\hess{\hbox{Hess}\,}

\def\e{\varepsilon}
\def\p{\partial}

\def\F{{\mathcal F}}

\begin{document}

\begin{frontmatter}

\title{Recovering the isometry type of a Riemannian manifold from local boundary diffraction travel times\tnoteref{thanks}}
  
\tnotetext[thanks]{This research was supported by National Science Foundation grant CMG DMS-1025318, the members of the Geo-Mathematical Imaging Group at Purdue University, the Finnish Centre of Excellence in Inverse Problems Research, Academy of Finland project COE 250215, the Research Council of Norway, and the VISTA project. The research was initialized at the Program on Inverse Problems and Applications at MSRI, Berkeley, during the Fall 2010.}
  
  \author[PU]{Maarten V. de Hoop}
  \ead{mdehoop@math.purdue.edu}
  
  \author[PU]{Sean F. Holman\corref{cor}}
  \ead{sfholman@math.purdue.edu}
 \cortext[cor]{Corresponding author.}
 
 \author[NORSAR]{Einar Iversen}
\ead{Einar.Iversen@norsar.com} 

\author[HEL]{Matti Lassas}
\ead{matti.lassas@helsinki.fi}

\author[NUST]{Bj{\o}rn Ursin}
\ead{bjorn.ursin@ntnu.no}

\address[PU]{Department of Mathematics, Purdue University, 150 N. University Street, West Lafayette IN 47907, USA.}

\address[NORSAR]{NORSAR, Gunnar Randers vei 15, P.O. Box 53, 2027 Kjeller, Norway}

\address[HEL]{Department of Mathematics and Statistics, Gustaf Hallstromin katu 2b, FI-00014 University of Helsinki, Helsinki, Finland.}

\address[NUST]{Department of Petroleum Engineering and Applied Geophysics, Norwegian University of Science and Technology, S.P. Andersensvei 15A, NO-7491 Trondheim, Norway.}

\begin{abstract}

\indent We analyze the inverse problem, originally
formulated by Dix~\cite{dix} in geophysics, of reconstructing the wave
speed inside a domain from boundary
measurements associated with the single scattering of seismic waves. 
We consider {\mltext a domain $\tilde M$ with a varying and possibly anisotropic wave speed which we model as a Riemannian metric $g$. For our data, we assume that $\tilde M$ contains a dense set of point scatterers and that in a subset $U\subset \tilde M$, modeling a region
containing measurement devices, we can measure the wave fronts of the single scattered waves diffracted from the point scatterers. The inverse problem we study is to recover
the metric $g$ in local coordinates anywhere on a set $M \subset \tilde M$ up to an isometry (i.e. we recover the isometry type of $M$). To do this we show that the shape operators related to wave fronts
produced by the point scatterers within $\tilde M$} satisfy a certain system of differential equations which may be solved along geodesics of the metric. In this way, assuming we know $g$ as well as the shape operator of the wave fronts in the region $U$, we may recover $g$ in certain coordinate systems (e.g. Riemannian normal coordinates centered at point scatterers). This generalizes the well-known geophysical
method of Dix to metrics which may depend on all spatial variables and be 
anisotropic. In particular, 
the  novelty of this solution lies in the fact that  it can
be used to reconstruct the metric  also in the presence of the caustics.\\

\noindent {\bf R\'esum\'e}\\

Nous analysons le probl\`eme inverse, \`a l'origine formul\'e par Dix~\cite{dix} en g\'eophysique, consistant \`a reconstruire la vitesse d'onde dans un domaine \`a partir des mesures aux fronti\`eres associ\'ees \`a
la dispersion simple des ondes sismiques. Nous consid\'erons un domaine $\tilde M$ avec une vitesse d'onde variable et \'eventuellement anisotrope que nous mod\'elisons par une m\'etrique Riemannienne $g$. Nous supposons que $\tilde M$ contient une densit\'e \'elev\'ee de points diffractants et que dans un sous-ensemble $U \subset \tilde M$, correspondant \`a un domaine contenant les instruments de mesure, nous pouvons mesurer les fronts d'onde de la diffusion simple des ondes diffract\'ees depuis les points diffractant. Le probl\`eme inverse que nous \'etudions consiste \`a reconstruire la m\'etrique $g$ en coordonn\'ees locales sur l'ensemble $M \subset \tilde M$ modulo une isom\'etrie (i.e. nous reconstruisons le type d'isom\'etrie). Pour ce faire nous montrons que l'op\'erateur de forme relatif aux fronts d'onde produits par les points diffractants dans $M$ satisfait un certain syst\`eme d'\'equations diff\'erentielles qui peut \^etre r\'esolu le long des g\'eod\'esiques de la m\'etrique. De cette mani\`ere, en supposant que nous connaissons $g$ ainsi que l'op\'erateur de forme des fronts d'onde dans la r\'egion $U$, nous pouvons retrouver $g$ dans un certain syst\`eme de coordonn\'ees (e.g. coordonn\'ees normales Riemannienne centr\'ees aux points diffractant). Ceci g\'en\'eralise la m\'ethode g\'eophysique de Dix \`a des m\'etriques qui peuvent d\'ependre de toutes les variables spatiales et \^etre anisotropes. En particulier, la nouveaut\'e de cette solution est de pouvoir \^etre utilis\'ee pour reconstruire la m\'etrique, m\^eme en pr\'esance des caustiques.
\end{abstract}



 
\end{frontmatter}
 
\pagestyle{myheadings}
\thispagestyle{plain}

\section{Introduction: Motivation of the problem}
\label{sec:1}


We consider a Riemannian manifold, $(M,g)$, of dimension $n$ with
boundary $\partial M$. We analyze the inverse problem, originally
formulated by Dix \cite{dix}, aimed at reconstructing $g$ from boundary
measurements associated with the single scattering of
seismic waves.  {\mltext When the waves produced by a source $F$ are modeled 
by the solution of the wave equation $(\p_t^2-\Delta_g)u(x,t)=F(x,t)$
on $(M,g)$, the geodesics $\gamma_{x,\eta}$ on $M$ correspond to the rays following the
propagation of singularities by the parametrix corresponding with the
wave operator on $(M,g)$ and the metric distance $d(x_1,x_2)$ of the points $x_1,x_2\in M$
corresponds to the travel time of the waves from the point $x_1$ to the point $x_2$. 
The phase velocity in this case is given by
$\mathrm{v}(x,\alpha) = [\sum_{j,k=1}^n g^{jk}(x) \alpha_j
  \alpha_k]^{1/2}$, with $\alpha$ denoting the phase or cotangent
direction. 

Below, we call the sets
 $\Sigma_{t,y}=\{\gamma_{y,v}(t);\ v\in T_yM,\
\|v\|_g=1\}$ generalized metric spheres (i.e. the images of the spheres $\{\xi\in T_yM;\ \|\xi\|_g=t\}$ in the tangent
space of radius $t$ under the exponential map). The mathematical formulation of Dix's problem is then the following: Assume
 that 
 \begin{enumerate} \item
 We are given an open set $\Gamma\subset \p M$, the metric tensor $g_{jk}|_{\Gamma}$ of the boundary, and normal derivatives
 $\p_\nu^p g_{jk}|_{\Gamma}$ for all $p\in \Z_+$, where $\nu$ is the normal vector of $\p M$
 and $g_{jk}$ is the metric tensor in the boundary normal coordinates.

\item 
For all $x\in \Gamma$ and $t>0$, we  are given at the point $x$ the second fundamental form of the generalized metric sphere  of $(M,g)$
 having the center $y_{x,t}$ and radius $t$.  Here,  $y_{x,t}=\gamma_{x,\nu(x)}(t)$ is the end point
 of the geodesic that starts from $x$ in the $g$-normal direction $\nu(x)$ to $\p M$ and has the length $t$.
   \end{enumerate}
 
 \noindent Dix's inverse problem is the question if one can use these data to determine uniquely
  the metric tensor $g$ on the set $W\subset M$ that can be connected
  to $\Gamma$ with a geodesic that does not intersect with a boundary.
 
 The above problem is the mathematical idealization of an imaging
problem encountered in geophysics
where the goal
is  to determine
the speed of waves (e.g.\ pressure waves) in a body from the external measurements. Roughly speaking, in 
Dix's inverse problem we assume that we observe  wave fronts of the waves
reflected from a large number of point scatterers, but for a given wave front
we do not know from which scatterer it reflected. Let us describe this
problem next in more detail:
 Assume that  the domain $M$ that contains a quite dense set of point scatterers
(called also diffraction points). Let
$\Gamma\subset \p M$ be the  acquisition surface on which
we have sources and measure the scattered waves. 
Consider a point source at the point  $x\in \Gamma$
that sends at time zero a wave that propagates into the domain and scatters from the point diffractors.
In the single scattering approximation, all these point scatterers can be considered as  new point sources.
We observe on $\Gamma$ the sum of the waves produced by these new point sources
and the observed wave fronts coincide with the generalized metric spheres of $(M,g)$.  
The wave fronts observed at the point $x$ and the time $2t$ are produced by the
point scatterers that can be connected to  $x$ with 
a geodesic whose length in the travel time metric is $t$.
Combining  the measurements corresponding to several point sources on $\Gamma$
to simulate a measurement that would be obtained using a wave packet
sent from $\Gamma$ in the normal direction one can determine for given $x\in \Gamma$ and $t$ the second fundamental form
 (or
equivalently, the shape operator) of the first wave front of the wave
produced by the point scatterer at $y_{x,t}$ that is observed
at $x$ to propagate
to the direction $-\nu(x)$. The wave speed in $W$ should be then reconstructed
using these shape operators.

Also, we note that in some cases when in $M$ there are  surface discontinuties
instead of point scatterers, one can use 
 signal processing of the measured data to produce the shape operators
 corresponding to point scatterers, see   \cite{GJI}.
This recovery of point scatterer data, which we will not describe in detail here, involves
 differentiating with respect to the offset between source and receiver locations in $\partial M$.
}

%
%
%

Earlier, Dix \cite{dix} developed a procedure, with a formula, for
reconstructing wave speed profiles in a half space  
$\R\times \R_+$ with an isotropic metric
that is one-dimensional (i.e. depends only on the depth coordinate). 
Despite this rather large restriction, Dix'
algorithm has played a crucial role in imaging problems in Earth sciences. 
We generalize this approach to the case of
multi-dimensional manifolds with general non-Euclidean metrics. Before continuing we add that since Dix, various adaptions have
been considered to admit more general wavespeed functions in a half
space. Some of these adaptions include the work of Shah \cite{shah}, Hubral \& Krey
\cite{hubralK}, Dubose, Jr. \cite{dubose}, Mann \cite{mann}, and Iversen
\& Tygell \cite{tygelI}. In the case that direct travel times are measured, rather than quantities related to point diffractors, the related mathematical formulations are either the boundary or lens rigidity problems. Much work has been done on these problems (see \cite{SU-lens,SU-review} and the references contained therein).  

\begin{figure}
\centering
\includegraphics[width=100mm]{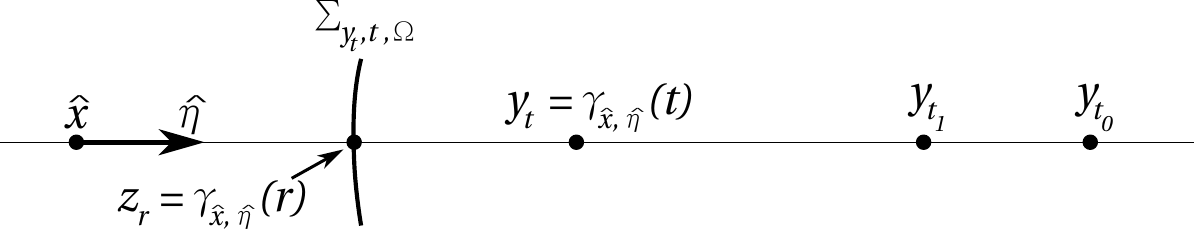}
\caption{Notation used throughout the paper, following the geodesic
         $\gamma$.}
\label{fig:notation}
\end{figure}

\section{Introduction: Definitions and main results}
\label{sec:def}

 \subsection{Background and notation}

Before continuing, we briefly mention some general references to Riemannian geometry \cite{Eisen,
  Petersen-book, Petersen-bulletin}. {\mltext As this paper
  is intended also for researchers working on applied sciences, we recall some standard notations and constructions}
in local coordinates, $(x^1,x^2,\dots,x^n)$. The metric tensor is
given by $g_{jk}(x) \dd x^j \dd x^k$ and the inverse of the matrix
$[g_{jk}]$ is denoted by $[g^{jk}]$. {\mltext Throughout the paper 
we use Einstein summation convention, summing over indexes that
appear both as sub- and super-indexes.} The Riemannian curvature tensor,
$R_{i j kl}$, is given in coordinates by
\[
   R^i_{j kl} = \frac{\p }{\p x^k } \Gamma^i _{j l}
        - \frac{\p }{\p x^l} \Gamma^i _{j k}
   + \Gamma^p_{j l} \Gamma^i_{p k} - \Gamma^p_{j k} \Gamma^i_{p l} ,
   \quad R^p_{j kl} = g^{pi} R_{ijkl} ,
\]
where $\Gamma^i_{j k}$ are the Christoffel symbols,
\[
   \Gamma^i _{j k} = \frac 12 g^{pi} \left(\mstrut{0.4cm}\right.
   \frac{\p g_{j p}}{\p x^k} + \frac{\p g_{k p}}{\p x^j}
       - \frac{\p g_{j k}}{\p x^p} \left.\mstrut{0.4cm}\right) .
\]
When $X,Y\in T_xM$ are vectors, then the curvature operator $R(X,Y):T_xM
\to T_xM$ is defined by the formula
\ba
g(R(X,Y)V,W)=R_{i j kl}X^iY^jV^kW^l,\quad V,W\in T_xM.
\ea
Finally, $\nabla_k = \nabla_{\p_k}$ is the covariant
derivative in the direction $\p_k = \frac\p{\p x^k}$, which is defined for
a (1,1)-tensor field $A^j_l$ by
\[
   \nabla_k A^j_l = \frac{\p}{\p x^k} A^j_l
       - \Gamma^p_{kl} A^j_p + \Gamma^j_{kp} A^p_l ,
\]
and for a (1,0)-tensor field $B^l$ and a (0,1)-tensor field $B_l$ by
\[
   \nabla_k B^l = \frac{\p}{\p x^k} B^l
       + \Gamma^l_{kp} B^p ,\quad
   \nabla_k B_l = \frac{\p}{\p x^k} B_l
       - \Gamma^p_{lk} B_p .
\]
If $f \in C^\infty$ then the gradient of $f$ with respect to $g$ is a (1,0)-tensor field (i.e. a vector field) given in coordinates by
\[
(\nabla f)^l = g^{lj} \frac{\partial f}{\partial x^j}.
\] 

\subsection{Description of the problem and results}

{\mltext Let us next formulate rigorously the setting for a modification of the above Dix's inverse problem that we will study in this paper.
Let $(M,g)$ be a $C^\infty$-smooth Riemannian manifold with boundary.}
We introduce an extension, $(\tilde M,\tilde g)$, $M \subset \tilde
M$, of $(M,g)$ which is a complete or closed manifold containing $M$ so that
$\tilde g|_M = g$. For simplicity we simply write $\tilde g = g$ and
assume that we are given $U=\tilde M \setminus M$ and the metric $g$ on
$\tilde M \setminus M$. We note that if  $M$ is compact
and the boundary of $M$ is convex,
the travel times between boundary points determine the normal
derivatives, of all orders, of the metric in boundary normal
coordinates \cite{LassasSU:2003}; hence, in the case of a convex
boundary, the smooth extension can be constructed when the travel
times between the boundary points are given. 

In the following we consider  a complete or closed Riemannian manifold $\tilde M$ with the
measurement data given in an open subset $U\subset \tilde M$.
{\mltext The tangent and cotangent bundles of $\tilde M$ at $x\in \tilde M$ are
denoted by $T_x\tilde M$ and $T_x^*\tilde M$, and
the unit vectors at $x$ are denoted by $\Omega_x\tilde M=\{v\in T_x\tilde M;
\ \|v\|_g=1\}$. We denote by $\exp_x:T_x\tilde M\to \tilde M$ the 
exponential map of $(\tilde M,g)$ and when $\eta\in \Omega_x\tilde M$,
we denote the geodesics  having the initial values $(x,\eta)$
by $\gamma_{x,\eta}(t)=\exp_x(tv)$.
Let $B_{\tilde M}(y,t)=\{x\in \tilde M;\
d(x,y)<t\}$ denote the metric ball of radius $t$ and center $y$ and call
the set $\{\exp_y(\xi)\in \tilde M;\ \|\xi\|_g=t\}$ the 
generalized metric sphere.}

We now express the \textit{data} described in the previous section in more precise terms. Eventually it will be related to the shape operators of so-called spherical surfaces.

\bigskip

\begin{definition} Let $U \subset \tilde M$ be an open set. The 
family ${\mathcal S}_U^{or}$ of oriented spherical surfaces is the set of all triples
$(t,\Sigma,\nu)$ satisfying the following properties
\begin{itemize}

\item [(i)] $t>0$,

\item  [(ii)] $\Sigma \subset U$ is a non-empty connected
$C^\infty$-smooth $(n-1)$-dimensional submanifold,

\item  [(iii)] There exists a $y \in \tilde M$ and an open set $\Omega \subset \Omega_y \tilde M$ such that
\begin{equation} \label{def: Sigma}
   \Sigma=\Sigma_{y,t,\Omega} = \{ \gamma_{y,\eta}(t) ;\ \eta \in \Omega \}.
\end{equation}

\item  [(iv)] $\nu$ is the unit normal vector field on $\Sigma$ given by
\beq \label{def: nu}
\nu(x)=\dot \gamma_{y,\eta}(t)\quad\hbox{at the point }
x= \gamma_{y,\eta}(t).
\eeq 
\end{itemize}
{\mltext The family ${\mathcal S}_U$ of spherical surfaces
is the set 
$${\mathcal S}_U=\{(t,\Sigma);\
\hbox{there exists }(t,\Sigma,\nu)\in {\mathcal S}_U^{or}\},
$$
that is, it contains the  same information as ${\mathcal S}_U^{or}$ but
not the orientation of the spherical surfaces.}
\end{definition}

\medskip
{\mltext
Note that if $(t,\Sigma,\nu)\in {\mathcal S}_U^{or}$ then in (\ref{def: Sigma})
the set $\Omega$ can be written in the form
$\Omega=\{\dot\gamma_{x,\nu(x)}(-t);\ x\in\Sigma\}$ and as $\Sigma$ is 
connected, thus also $\Omega$ is connected.}
\medskip

\noindent
{\mtext In reflection seismology one refers to $\Sigma$ as the (partial)
front of a point diffractor. Note that in the definition of ${\mathcal S}_U$
we consider arbitrary $y\in \tilde M$, including points $y$ in $U$.
Due to this, we have the following result stating that
${\mathcal S}_U$ determines both the metric $g$ in $U$ and
the wave front data with orientation, that is, ${\mathcal S}_U^{or}$.
Even though the determination of the metric in $U$ is not very interesting
from the point of view of applications, we state the proposition for
mathematical completeness.

 \medskip

\begin{proposition}\label{prop: determination in U}
Assume that we are given the open set $U$ as a differentiable manifold and
the family of spherical surfaces
${\mathcal S}_U$. These data determine
the metric $g$ in $U$ and the family of the   oriented spherical surfaces
${\mathcal S}_U^{or}$.
\end{proposition}
 \medskip

Proposition \ref{prop: determination in U} is proven in the Appendix.


Let ${x} \in U$ and ${\eta} \in \Omega_{{x}}\tilde M$.
We say that $(t,\Sigma,\nu)\in {\mathcal S}_U^{or}$ is associated
to  the pair $(x, \eta)$ if $x\in \Sigma$ and $\nu(x)=-\eta$.
It is easy to see that if $(t,\Sigma,\nu)\in {\mathcal S}_U^{or}$ is associated
to  the pair $(x, \eta)$ then 
we can represent $\Sigma$ in the form (\ref{def: Sigma}) where
 $y=\gamma_{x,\eta}(t)$ and $\Omega \subset \Omega_y \tilde M$
 is such that ${\zeta} = - \dot \gamma_{x,\eta}(t)\in \Omega$.

%
%

Once again, suppose that ${x}\in U$, ${\eta} \in \Omega_{{x}} \tilde M$.
Now we proceed with more geometrical constructions along  the geodesics $\gamma_{x,\eta}$. Let $F_k(r)=F_k(x,\eta,r)$, $k=1,2,\ldots,n$ be a linearly
independent and parallel set of vector fields defined on $\gamma_{x,\eta}(\R)$.
This means that $F_k(x,\eta,r)\in T_{\gamma_{x,\eta}(r)}\tilde M$ and $\nabla_{\dot\gamma_{x,\eta}(r) }F_k(r) = 0$. 
We assume that $F_n(x,\eta,r)= \dot \gamma_{x,\eta}(r)$.
Denote by $\hat g_{jk}$ the inner products
\begin{equation}\label{hat g matrix}
   \hat g_{jk} = g(F_j(r),F_k(r)) .
\end{equation}
Because the vector fields $F_j$ are parallel, $\hat{g}_{jk}$ does not depend on $r$. 
Let $f^j=f^j(x,\eta,r)$, $j=1,\ldots,n$ be the co-frame dual to $F_j$. This means that
\[
   \bra f^j(r),F_k(r) \cet = \delta_k^j
\]
where $\langle .,. \rangle$ denotes the usual pairing of
  $T_y \tilde M$ and $T_y^* \tilde M$. Let $\Psi=\Psi_{x,\eta}:\R^n\to \tilde M$ be the map
\ba
\Psi_{x,\eta}(s^1,s^2,\dots,s^{n-1},r)=\exp_{q(r)} \left ( \sum_{k=1}^{n-1}s^kF_k({x,\eta},r) \right ), \quad
\hbox{where }q(r)=\gamma_{x,\eta}(r).
\ea
For all $r\in \R$ the point $q(r)$ has a neighborhood $B_{\tilde M}(q(r),\e)$, $\e>0$
so that there exists an smooth inverse map $\Psi_{x,\eta}^{-1}:B_{\tilde M}(q(r),\e)\to \R^n$.
We call such inverse maps the Fermi coordinates.

Let $R$ be the curvature operator of $(\tilde M,g)$. Below, we denote 
\beq\label{eq: r coefficients}
\mathbf{r}_j^k(x,\eta,r) = \langle f^k(x,\eta,r), R(F_j(x,\eta,r), \dot \gamma_{x,\eta}(r) ) \dot \gamma_{x,\eta}(r)  \rangle
\eeq
and call ${\bf r}^k_j(x,\eta,r)$ the {\it curvature coefficients} of the frame $(F_j(x,\eta,r))_{j=1}^n$.

For $t\geq 0$, let  ${\mathcal C}(x,\eta,t)$ be the set of those $r\geq 0$ for
which $z_r=\gamma_{x,\eta}(r)$ and $y_t=\gamma_{x,\eta}(t)$ are conjugate points
on the geodesic $\gamma_{x,\eta}$.
When $t > r \geq 0$ and $r \not\in {\mathcal C}(x,\eta,t)$,
there is a spherical surface $\Sigma_{r,t}:=\Sigma_{t-r,y_t,\Omega}$
of the form (\ref{def: Sigma}) with some neighborhood  $\Omega\subset \Omega_{y_t}\tilde M$
of $\zeta=- \dot\gamma_{x,\eta}(t)$.
%
Then $z_r=\gamma_{x,\eta}(r)\in \Sigma_{r,t}$. 
In this case we write the shape operator of $\Sigma_{r,t}$ at $z_r$ as $S_{x,\eta,r,t}$. Thus $S_{x,\eta,r,t} \in (T_1^1)_{z_r} \tilde M$ is defined by 
\[
S_{x,\eta,r,t} X = \nabla_X \nu
\]
for all $X\in T_{z_r}\tilde M$ where $\nu$ is the normal vector field for $\Sigma_{r,t}$ satisfying (\ref{def: nu}). Let us write the shape operator with respect to the parallel frame:
\beq\label{S-coefficients}
S_{x,\eta,r,t} = {\bf s}^k_j(x,\eta,r,t) f^j(x,\eta,r) \otimes F_k(x,\eta,r).
\eeq
We say that ${\bf s}^k_j(x,\eta,r,t)$ are
{\it the coefficients of the second fundamental forms of the spherical
surfaces} on the geodesic $\gamma_{x,\eta}$ corresponding
to the frame $F_k(0)$. 
{\mltext 
The family of the oriented spherical surfaces
${\mathcal S}_U^{or}$ and the metric tensor $g$ in $U$
(which are in fact determined by ${\mathcal S}_U$ by Proposition \ref{prop: determination in U})),
determine all 
triples  $(t,\Sigma,\nu)\in {\mathcal S}_U^{or}$ that are
associated to  the pair $(x,\eta)$. Note that if $t>0$ is such
that   $x$ and $\gamma_{x,\eta}(t)$ are not conjugate along $\gamma_{x,\eta}$,  there exists
at least one triple $(t,\Sigma,\nu) \in \mathcal{S}_U^{or}$ that is associated to $(x,\eta)$.
Note that all such surface $\Sigma$ have the 
same the shape operator $S_{x,\eta,r,t}$ with $r=0$   at $x$. 
Thus, by computing  the shape operator  of such surfaces $\Sigma$ at $x$ we can find 
the operator $S_{x,\eta,0,t}$ and furthermore
the coefficients
${\bf s}^k_j(x,\eta,0,t)$ for all $t>0$ such that $\gamma_{x,\eta}(t)$ is not conjugate point to $x$.}


Our main result is the following.

\medskip\medskip

\begin{theorem} \label{thm: main}
Let $(\tilde M,g)$ be a complete or closed Riemannian manifold of dimension $n$
and $U\subset \tilde M$ be open. Then \medskip 

(i)  Let $x\in U$, $\eta\in \Omega_x\tilde M$ and $F_k(r)$, $k=1,2,\dots,n$
be linearly independent parallel vector fields along $\gamma_{x,\xi}$.   
Assume that we are given $\hat g_{jk}=g(F_j,F_k)$ for $j,k=1,2,\dots,n$ and 
\mltext the coefficients of the second fundamental forms of the spherical
surfaces corresponding to the frame $F_k(0)$, that is, ${\bf s}^k_j(x,\eta,0,t)$, for all $t\in \R_+\setminus \mathcal C(x,\eta,0)$. Then we can determine uniquely ${\bf s}^k_j(x,\eta ,r,t)$ for all 
 $t > 0$ and $r<t$, $r\in \R_+\setminus \mathcal C(x,\eta,t)$,
and the curvature coefficients $\mathbf{r}_j^k(x,\eta ,r)$ for all $r\in \R_+$.

\medskip

\noindent
Consequently,

\medskip 

(ii) Assume that we are given the open set $U \subset \tilde M$, metric $g$ on $U$,
$x_0\in U$ and a unit vector $\eta_0\in \Omega_{x_0}\tilde M$ and let
$\mathcal V$ be a neighborhood of $(x_0,\eta_0)$ in $T\tilde M$. Moreover,
{\mnewtext  assume
we are given  the set 
\ba
 {\mathcal S}_{U,\V}^{or}:=\{(t,\Sigma,\nu)\in {\mathcal S}_U^{or};\ (x,-\nu(x))\in \V\hbox{ for
 all }x\in \Sigma\}.\ea 
Then for all $r>0$ there is $\rho=\rho(x_0,\eta_0,r)$ such that the Fermi
coordinates $\Psi_{x_0,\eta_0}^{-1}$ 
associated to the geodesic $\gamma_{x_0,\eta_0}(\R)$
are well defined in an open set
\[
V_{x_0,\eta_0,r}=
\Psi_{x_0,\eta_0}(B_{\R^{n-1}}(0,\rho) \times  
(r-\rho,r+\rho)),
\]
and the above data
determine uniquely $(\Psi_{x_0,\eta_0})_*g$, that is, the metric $g$ in Fermi coordinates in $
V_{x_0,\eta_0,r}$. This is the meaning of ``reconstruction of the isometry type of the metric in the Fermi coordinates
near $\gamma_{x_0,\eta_0}(\R_+)$."}

\end{theorem}

\medskip\medskip

In the above theorem, (i) says that the shape operators $S_{x,\eta,r,t}$ of the $\gamma_{x,\eta}(t)$
centered generalized spheres with radius $t-r$ can be uniquely determined from the shape operators  $S_{x,\eta,0,t}$
corresponding to $r=0$ when the metric in $U$ is known. The claim (ii)
says that the Riemannian metric near the geodesic $\gamma_{x_0,\eta_0}(\R_+)$
can be determined from the knowledge of wave fronts propagating
close to this geodesic. From the point of view of applications, it is particularly
important that the reconstruction can be done past the first conjugate point,
that is, beyond the caustics of the reflected waves.


\medskip

We point out that the reconstruction method we develop in this paper
is constructive and that it is based on solving a system of ordinary
differential equations which are satisfied by $\p_t^p{\bf s}^k_j(x,\eta ,r,t)$, $p=0,1,2,3$ along each geodesic. 

{\mltext 
Using Theorem \ref{thm: main} we can prove the unique determination
of the universal covering space.

\medskip

\begin{theorem} \label{thm: main 2}
Let $\tilde M$ and $\tilde M^\prime $ be two smooth  (compact or complete)
Riemannian manifolds and  $U\subset \tilde M$ and $U^\prime \subset 
\tilde M^\prime $  be such non-empty open sets  that
there is a diffeomorphism $\Phi:U \to U^\prime $ and
${\mathcal S}_{U^\prime }=\{(t,\Phi(\Sigma));\  (t,\Sigma)\in {\mathcal S}_U\}$. 
Then there is a Riemannian manifold $(N,g_N)$ such that there are
 Riemannian covering maps $F:N\to \tilde M$ and   $F^\prime :N\to \tilde M^\prime $, that is,
$\tilde M$ and $\tilde M^\prime $ have isometric
universal covering spaces.
\end{theorem}
\medskip

The next example shows that ${\mathcal S}_U$ does not determine
the manifold $(\tilde M,g)$ but only its universal covering space.
\medskip

 {\bf Example 1.} Let $(\tilde M,g)$ be the flat torus $\R^2/\Z^2$,
 and $(\tilde M',g)$ be the flat torus $\R^2/(2\Z\times 2\Z)$.
  Below, we consider  $(\tilde M,g)$ and
 $(\tilde M^\prime ,g^\prime )$ as the squares $[0,1]^2$ and $[0,2]^2$
 with the parallel sides being glued together. 
 Both $(\tilde M,g)$ and $(\tilde M^\prime ,g^\prime )$ have the universal covering
 space $\R^2$ with the Euclidean metric.
 Let $U$ and $U^\prime $ be the 
 disc $B(p,\frac 14)$ of the radius $\frac 14$ and the center $p=(\frac 12,\frac 12)$.
 We see that  both the collection ${\mathcal S}_U^{or}$ for $(\tilde M,g)$ and the collection ${\mathcal S}_{U^\prime }^{or}$ for $(\tilde M^\prime ,g^\prime )$ consist of
 triples $(\Sigma,t,\nu)$ where $t>0$, $\Sigma$ is a connected circular arc having radius $t$  
 that is a subset of the disc $B(p,\frac 14)$, that is, $\Sigma=\{(t\sin \alpha + x,t\cos\alpha+y)\in U;\ |\alpha-\alpha_0|<c_0\}$, and $\nu$ is the exterior normal vector of $\Sigma$.
Note that the circular arc $\Sigma$ can also be a whole circle if $t$ is small enough.
This shows that the knowledge of $U$ and ${\mathcal S}_{U^\prime }^{or}$ is not
always enough to determine uniquely the manifold $(\tilde M,g)$ but only its universal covering space. 
\medskip

The above example shows that
 ${\mathcal S}_{U^\prime }^{or}$ contains less information than many other data
sets used in inverse problems that determine the manifold
$(\tilde M,g)$ uniquely such as the hyperbolic Dirichlet-to-Neumann map
$\Lambda^{wave}$, the parabolic Dirichlet-to-Neumann map 
associated to the heat kernel
considered in \cite{BK,KK,KKL,KKLM,Oksanen,StU1,StU2}, the  
boundary distance representation $R(M)$ (i.e., the boundary distance
functions) considered in \cite{AKKLT,Ku,KKL}, or
the broken scattering relation considered in \cite{KLU}; see also related
data sets in \cite{Uhl}.

 }

\subsection{Jacobi and Riccati equations}
\label{sec:2}

Before moving to the actual reconstruction procedure we collect a few more geometrical formulae that will be useful. {\mltext First, if we fix the initial 
data $(x,\eta)$ for the geodesic and a $t>0$, then $ S_t(r)=S_{x,\eta,r,t} $ can be thought of as a $(1,1)$-tensor field on the geodesic $\gamma_{x,\eta}$.  Let
$y_t=\gamma_{x,\eta}(t)$, $z_r=\gamma_{x,\eta}(r)$,
and $\zeta=-(t-r)\dot\gamma_{x,\eta}(t)$
so that $z_r=\exp_{y_t}(\zeta)$.
Assume that $y_t$ and $z_r$
 are not conjugate
points along $\gamma_{x,\eta}(t)$. Then the exponential function $\exp_{y_t}:T_{y_t}\tilde M\to \tilde M$ has
a local inverse $F_t=\exp_{y_t}^{-1}$ in a neighborhood $V$ of $z_r$
and the function $f_t(z)=\|F_t(z)\|_g$ is a generalized distance
function (i.e. $\|\nabla f_t(z)\|_g=1$, $z\in V$). As
the spherical surface
$\Sigma_{y,t-r,\Omega}$ near $z_r$
can be written as a level set of a generalized distance function $f_t$, it follows from the radial curvature (Riccati) equation \cite[sect. 4.2, Thm 2]{Petersen-book} that
\begin{equation} \label{eq:Sy1}
   -\nabla_{\p_r} S_t(r) + S_t(r)^2 = -R_{\p_r}(r),
\end{equation}
where $R_{\p_r}(r):T_{z_r}\tilde M\to T_{z_r}\tilde M$,
$R_{\p_r}(r):V \mapsto R(V,\p_r) \p_r$ is the so-called directional curvature operator associated with the Riemannian
curvature $R$ of $(\tilde M,g)$ and $\p_r=-\nabla f_t$. Note that 
at  the point $\gamma_{x,\eta}(r)$ we have 
$\p_r=\dot\gamma_{x,\eta}(r)$.}   Let 
$\mathbf{r}_j^k(r)=\mathbf{r}_j^k(x,\eta,r)$ be the coefficients defined in 
 (\ref{eq: r coefficients}), that is, 
\[
\mathbf{r}_j^k(r) = \langle f^k, R(F_j, \p_r) \p_r \rangle.
\] 
If we also express $S_t(r)=S_{x,\eta,r,t}$ in the parallel frame as in (\ref{S-coefficients})
with $\mathbf{s}_j^k(r,t)=\mathbf{s}_j^k(x,\eta,r,t)$, then equation (\ref{eq:Sy1}) becomes
\begin{equation} \label{Matrix Riccati}
- \partial_r \mathbf{s}_j^k + \mathbf{s}^k_p \mathbf{s}_j^p = - \mathbf{r}_j^k.
\end{equation}
Also there is an equation we will use relating Jacobi fields along $\gamma_{x,\eta}$ to $S_t$. A Jacobi field $J(r)$ along the geodesic
$\gamma_{x,\eta}$ is a vector field satisfying
\begin{equation} \label{eq:Jeq}
   \nabla_{\p_r}^2 J + R(J,\p_r) \p_r = 0.
\end{equation}
Writing
\[
J = J^k F_k
\]
this equation is
\[
\partial_r^2 J^k + \mathbf{r}^k_j J^j = 0.
\]
Finally, it follows from \cite[p. 36]{Gromoll} that if $J |_{r=t} = 0$,
then $\nabla_{-\p_r} J = S_t\,J$, that is,
\begin{equation} \label{eq:Jy1}
   -\nabla_{\p_r} J = S_t \, J.
\end{equation}
With respect to the parallel frame equation (\ref{eq:Jy1}) reads as
\begin{equation} \label{Matrix shape}
-\p_rJ^k = \mathbf{s}^k_j J^j.
\end{equation}

Our strategy for reconstruction is to show that from the data we can
reconstruct shape operator $S_{x,\eta,r,t}$ along each $\gamma_{x,\eta}$ using the Riccati equation
(\ref{eq:Sy1}). Then the Jacobi fields may be calculated using (\ref{eq:Jeq}), and since the Jacobi fields are the coordinate vectors for the local coordinates introduced in the next section this then allows recovery of the metric with respect to those coordinates.

\section{Coordinates associated to a spherical surface $\Sigma_0$}
\label{sec:3}

{\mltext
We now introduce a set of coordinates in which we will initially recover the metric $g$. Suppose that we fix $x_0 \in  U=\tilde M \setminus M$ and $\eta_0 \in T_{x_0} \tilde M$, and then pick a large $t_0 >0$ 
such that $x_0$ and $\gamma_{x_0,\eta_0}(t_0)$ are not conjugate points along $\gamma_{x_0,\eta_0}$. We will use the notation $y_{t_0} = \gamma_{x_0,\eta_0}(t_0)$. Suppose that $\Sigma_0\subset U$ is a spherical surface containing $x_0$ given by $y_{t_0}$ and $t_0$ in the form (\ref{def: Sigma}). When $\nu=\nu(x)$ is 
the normal vector field of $\Sigma_0$ oriented so that $\nu(x_0)=-\eta_0$,
we have 
$(\Sigma_0,t_0,\nu)\in {\mathcal S}_U^{or}$

Let us take arbitrary coordinates $\hat X:\hat V\to \R^{n-1}$ on $\hat V\subset \Sigma_0$
such that $\hat V$ is a neighborhood of $x_0$ and 
$\hat X(x_0)=0$. To slightly simplify the notation, we identify $\hat V$ with its image
$\hat X(\hat V)\subset \R^{n-1}$ and
the points $x\in \hat V$ with their coordinates 
$\hat{x} = (\hat{x}^1,\ldots,\hat{x}^{n-1})=\hat X(x)$. Also, denote
the geodesic starting normally to $\Sigma$ from $\hat X^{-1}(\hat x)\in \Sigma$ by $\gamma_{\hat{x}}(t) = \gamma_{\hat{x},-\nu(\hat{x})}(t)$. Let
$\iota:\Sigma_0\to \tilde M$ be the identical embedding. 
For $\hat x\in \hat X( \hat V)$, we define at the point $q=\hat X^{-1}(\hat x)$ the vectors
\ba
F_j(\hat x)=\iota_*\left (\frac \p{\p \hat X^j}\right ),\quad j\leq n-1,\quad \hbox{and}\quad
F_n(\hat x)=-\nu(\hat X^{-1}(\hat x))
\ea
Then $(F_j(\hat x))_{j=1}^n$ is a basis for $T_q\tilde M$.
Let $F_j(\hat x,r)$ be the parallel translation
of the vector $F_j(\hat x)$ along the geodesic $\gamma_{\hat x}(r)$.

Let  ${\mathcal C}(\hat x)$ be the set of those $r\in \R_+$ for
which $\hat x\in \Sigma$ and $\gamma_{\hat x}(r)$ are conjugate points
on the geodesic $\gamma_{\hat x}$, and let
\ba
& &\mathcal W=\{(\hat x,r)\in \hat V\times \R_+;
\ \hat x\in \hat V,\ t\in  \R_+\setminus {\mathcal C}(\hat x)\},\\
& &W=\{\gamma_{\hat x}(r)\in \tilde M;\ (\hat x,r)\in \mathcal W\}.
\ea
Then the map 
\beq\label{eq: the map Xhat V}
X_{\hat V}:\mathcal W\to W\subset  \tilde M,\quad  X_{\hat V}(\hat{x},r)=\gamma_{\hat{x}}(r) 
\eeq
is a local diffeomorphism. Below, we use the local inverse maps of $X_{\hat V}$
as local coordinaes on $W$. The $(\hat{x},r)$
coordinates, basically, are Riemannian normal coordinates centered at
$y_{t_0}$, but parametrized in a particular way: $\hat{x}$ can be
thought of as a parametrization of part of the sphere of
radius $t_0$ in $T_{y_{t_0}} \tilde{M}$, and then $r$ corresponds to
the radial variable in $T_{y_{t_0}} \tilde{M}$. Note also that the coordinate vectors in these coordinates are Jacobi fields along the geodesics $\gamma_{\hat{x}}$.

Below, we also denote
\ba
& &{\bf s}^k_j(\hat x,r,t)={\bf s}^k_j(\hat X^{-1}(\hat x),\nu(\hat X^{-1}(\hat x)),r,t),\\
& &\mathbf{r}_j^k(\hat{x},r)=\mathbf{r}_j^k(\hat X^{-1}(\hat x),\nu(\hat X^{-1}(\hat x)),r).
\ea
}


\section{Reconstruction of the shape operator along one geodesic}
\label{sec: 3b}

{\mltext Next, we fix $\hat{x}$ and aim to reconstruct 
the shape operator along $\gamma_{\hat{x}}$. Throughout the section we will suppress the dependence of all quantities on $\hat{x}$ because it is considered to be fixed. 

To simplify notations, in this section we use the conventions $\gamma = \gamma_{\hat{x}}$
 and
\ba
& &{\bf s}^k_j(r,t)={\bf s}^k_j(\hat x,r,t)\\
& &\mathbf{r}_j^k(r)=\mathbf{r}_j^k(\hat{x},r).\ea
}

\subsection{A lemma concerning curvature}

Our first step toward the reconstruction will be to prove a lemma relating, roughly speaking, the inverse of the shape operator of spherical surfaces as their radius goes to zero with the directional curvature operator at the center of the spherical surfaces.

\bigskip

\begin{lemma} \label{lem:1}
Let ${\bf S}(r,t) = ({\bf s}^j_k(r,t))_{j,k=1}^{n-1}$ be given by the matrices
defined in (\ref{S-coefficients}), $t_1 > 0$ and $i_1$ be the injectivity
radius of $(\tilde M,g)$ at $\gamma(t_1)$. Let $t,r \in
[t_1-i_1/2,t_1]$ with $t>r$ and ${\bf K}(r,t) = {\bf S}(r,t)^{-1}$. Then
\begin{equation}\label{K - R equation}
   {\bf K}(r,t) = (t-r) \, I + \frac{(t-r)^3}{3} \, {\bf R}(t)
                            + \mathcal{O}((t-r)^4) ,\quad
   {\bf R}(t) = ({\bf r}^j_k(t))_{j,k=1}^{n-1} ,
\end{equation}
where $\mathcal{O}((t-r)^4)$ is estimated in a norm on the space of
matrices $\R^{n \times n}$.
\end{lemma}

\bigskip

\begin{proof}
Throughout this proof all the indices and sums will run from $1$ to $n-1$ unless otherwise noted. We use Jacobi fields ${\bf j}^k_{(m)}(s;r,t) F_k(s)$ on
$\gamma([r,t])$ satisfying
\begin{equation} \label{Matrix Jacobi 2}
   \p_s^2 {\bf j}^k_{(m)}(s;r,t)
          + {\bf r}^k_p(s ) {\bf j}^p_{(m)}(s ;r,t) = 0 ,\quad
   s  \in [r,t] ,
\end{equation}
supplemented with the boundary data
\begin{equation} \label{Matrix ptnsrc 2}
   {\bf j}^k_{(m)}(s ;r,t)|_{s =r} = \delta^k_m ,\quad
     {\bf j}^k_{(m)}(s ;r,t)|_{s =t} = 0 .
\end{equation}
We will consider these when $t-r = \e > 0$ is sufficiently
small. Next, we freeze $t$ and introduce notations
\begin{eqnarray}
   v^k_{(m)}(z,t) &=& {\bf j}^k_{(m)}(t-z;t-\e,t) ,
\label{eq:vkmdef}
\\
    \rho^j_k(z,t) &=& {\bf r}^j_k(t-z) ,
\end{eqnarray}
where $z \in [0,\e]$. Equations (\ref{Matrix Jacobi 2})-(\ref{Matrix
  ptnsrc 2}) attain the form
\begin{eqnarray}
   \p_z^2 v^k_{(m)}(z,t) + \rho^k_p(z,t) v^p_{(m)}(z,t)
          &=& 0 ,\quad z \in [0,\e] ,
\\
\nonumber
   v^k_{(m)}(z,t)|_{z=0} &=& 0 ,\quad
   v^k_{(m)}(z,t)|_{z=\e} = \delta^k_m .
\end{eqnarray}
We then let
\begin{eqnarray} \label{eq:wedef}
   w^k_{(m);\e}(y,t) &=& v^k_{(m)}(\e y,t) ,
\\
   \sigma^j_{k;\e}(y,t) &=& \rho^j_k(\e y,t) ,
\end{eqnarray}
where $y \in [0,1]$. We drop the subscripts and superscripts for
simplicity of notation and view $w$ and $\sigma$ as matrices. Then
\begin{eqnarray}\label{Jacobi for w}
   \p_y^2 w(y,t) + \e^2 \sigma(y,t) w(y,t)
          &=& 0 ,\quad y \in [0,1] ,
\\
\nonumber
   w(y,t)|_{y=0} &=& 0 ,\quad
   w(y,t)|_{y=1} = I .
\end{eqnarray}
The supremum of the norm of the Riemannian curvature tensor is
bounded in compact sets and thus we see that 
$\|\sigma(y,t)\|$ is uniformly bounded 
over $y \in [0,1]$ and $t \in [t_0-i_0/2,t_0]$. Thus we see that there are $C_1$ and 
$\e_1>0$ such that if $0<\e<\e_1$, then
\begin{equation} \label{c6 estimate}
  \| (\p_y^2 + \e^2 \sigma(y,t))^{-1} \|_{L^2([0,1])
           \to H^1_0([0,1])} \le C_1 ,
\end{equation}
for some constant $C_1 > 0$. This implies that there is a $C_2 > 0$
such that for all $t$ and $\e$ corresponding with $r,t \in
[t_0-i_0/2,t_0]$,
\begin{equation} \label{eq:wbound}
   \| w (\cdotp,t)\|_{H^1([0,1])} \le C_2 .
\end{equation}

We expand $\sigma(y,t)$,
\begin{equation} \label{eq:sigexp}
   \sigma(y,t) = \sigma(0,t) + \mathcal{E}_{\sigma}(y,t) ,
\end{equation}
with 
\begin{equation} \label{extra epsilon estimate}
   \| \mathcal{E}_{\sigma} (\cdotp,t) \|_{L^{\infty}([0,1])}
                 \le C_3 \e ,
\end{equation}
where $C_3$ depends on  the
supremum of $\|\nabla R\|_g$ and the norms of $f^k$ and $F_k$
on the geodesic $\gamma([0,t_0])$, i.e., the $C^3$-norm of  the metric. We expand $w$ accordingly,
\begin{equation}
   w(y,t) = w^0(y,t) + \mathcal{E}_w(y,t) ,
\end{equation}
where
\begin{eqnarray}
   \p_y^2 w^0(y,t) + \e^2 \sigma(0,t) w^0(y,t)
          &=& 0 ,\quad y \in [0,1] ,
\\
\nonumber
   w^0(y,t)|_{y=0} &=& 0 ,\quad
   w^0(y,t)|_{y=1} = I .
\end{eqnarray}
We observe that there is a constant $C_4 > 0$, such that for all $y \in [0,1]$ and $t
\in [t_1-i_1/2,t_1]$,
\begin{equation} \label{c4 estimate}
  \| (\p_y^2 + \e^2 \sigma(0,t))^{-1} \|_{L^2([0,1])
           \to H^1_0([0,1])} \le C_4 .
\end{equation}
This implies that there is a $C_5 > 0$ for all $t$ and $\e$ such that
$r,t \in [t_1-i_1/2,t_1]$,
\[
   \| w^0 (\cdotp,t)\|_{H^1([0,1])} \le C_5 .
\]
Now
\begin{eqnarray*}
   \p_y^2 \mathcal{E}_w(y,t)  + \e^2 \sigma(0,t) \mathcal{E}_w(y,t)
   &=& -\e^2\mathcal{E}_{\sigma}(y,t) w(y,t)  ,\quad y \in [0,1] ,
\\
\nonumber
   \mathcal{E}_w(y,t)|_{y=0} &=& 0 ,\quad
   \mathcal{E}_w(y,t)|_{y=1} = 0 .
\end{eqnarray*}
Using (\ref{extra epsilon estimate}), (\ref{eq:wbound}) and (\ref{c4
  estimate}), we find that there is a $C_6 > 0$ such that
\[
   \| \mathcal{E}_w (\cdotp,t) \|_{H_0^1([0,1])} \le C_6 \e^3
\]
for $\e$ sufficiently small.

Denoting $\lambda = \lambda(t) = \sqrt{\sigma(0,t)} \in \C^{(n-1)\times
  (n-1)}$ (we can use any branch of the matrix square root) we get
\[
   w^0(y,t) = [\sin(\e \lambda(t))]^{-1}
                        \sin(\e \lambda(t) y) .
\]
Expanding this solution in $\e$ yields
\[
   w^0(y,t) = y \, \left(\mstrut{0.4cm}\right. I
           + \frac 16 \e^2 \lambda(t)^2 (1 - y^2)
             + \mathcal{O}(\e^4) \left.\mstrut{0.4cm}\right)
\]
whence
\begin{equation} \label{eq:wexpe}
   w(y,t) = y \, \left(\mstrut{0.4cm}\right. I
           + \frac 16 \e^2 \lambda(t)^2 (1 - y^2)
             \left.\mstrut{0.4cm}\right) + \mathcal{E}_{w;1}(y,t) ,
\end{equation}
where $\| \mathcal{E}_{w;1} (\cdotp,t) \|_{H^1([0,1])} \le C_6^\prime  \e^3$,
and
\begin{equation} \label{eq:dwexpe}
   \p_y w(y) = I + \frac 16 \e^2 \lambda(t)^2
         - \frac 12 \e^2 \lambda(t)^2 y^2
                   + \p_y \mathcal{E}_{w;1}(y,t) ,
\end{equation}
where $\| \p_y \mathcal{E}_{w;1} (\cdotp,t) \|_{L^2([0,1])} \le C_7
\e^3$.

We recall (\ref{eq:vkmdef}) and differentiate,
\[
   \p_z v^k_{(m)}(z,t)|_{z=\e}
         = -\p_s {\bf j}^k_{(m)}(s;t-\e,t)|_{s=t-\e} .
\]
Because $v^k_{(m)}(z,t) |_{z=0} = {\bf j}^k_{(m)}(t;t-\e,t)= 0$,
\[
   -{\p_s} {\bf j}^k_{(m)}(s;t-\e,t)
       = {\bf s}^k_p(s,t) {\bf j}^p_{(m)}(s;t-\e,t)
\]
(cf.~(\ref{Matrix shape})). Using this identity at $s = t-\e$, we find
that
\begin{equation}
   \p_z v(z,t)|_{z=\e} = {\bf S}(t-\e,t) v(z,t)|_{z=\e} .
\end{equation}
However, $v(z,t)|_{z=\e} = I$, so that, with $\tilde{\bf S}(t -\e,t)
:= ( \hat{g}_{lk} {\bf s}^k_m(t -\e,t) )_{l,m=1}^{n-1}$,
\begin{eqnarray*}
   \tilde{\bf S}(t -\e,t) &=& {\bf S}(t -\e,t) v(\e,t) \, \cdotp v(\e,t)
   = \p_z v(z,t)|_{z=\e} \, \cdotp v(\e,t)
\\
   &=& \, \int_0^\e (\p_z v(z,t) \, \cdotp \p_z v(z,t)
           + \p_z^2 v(z,t) \, \cdotp v(z,t)) \, dz ,
\end{eqnarray*}
where $v \cdotp v$ stands for $v^k_{(l)} v^j_{(m)} \hat{g}_{kj}$.
Substituting (\ref{eq:wedef}) yields
\begin{eqnarray*}
   \tilde{\bf S}(t -\e,t) &=&
   \int_0^1 (\e^{-1} \p_y w(y,t) \, \cdotp \e^{-1} \p_y w(y,t)
      + \e^{-2} \p_y^2 w(y,t) \, \cdotp w(y,t)) \, \e dy
\\
   &=& \e^{-1} \int_0^1 (\p_y w(y,t) \, \cdotp \p_y w(y,t)
      - \e^2 \sigma(y,t) w(y,t) \, \cdotp w(y,t)) \, dy ,
\end{eqnarray*}
using Jacobi equation (\ref{Jacobi for w}).

Inserting expansions (\ref{eq:wexpe}) and (\ref{eq:sigexp}) gives
\begin{multline*}
   \tilde{\bf S}(t-\e,t) = \e^{-1} \int_0^1 \bigg(
   \bigg(I +
   \frac 16 \e^2 \lambda(t)^2 - \frac 12 \e^2 \lambda(t)^2 y^2
                       \bigg) \cdotp
   \bigg(I +
   \frac 16 \e^2 \lambda(t)^2 - \frac 12 \e^2 \lambda(t)^2 y^2
                       \bigg)
\\
   - \e^2 \sigma(0,t) y^2 \bigg(I +
   \frac 16 \e^2 \lambda(t)^2 (I-y^2) \bigg) \cdotp
   \bigg(I + \frac 16 \e^2 \lambda(t)^2 (I-y^2) \bigg)
   \bigg) dy + \mathcal{O}(\e^2)
\end{multline*}
so that
\begin{multline*}
   {\bf S}(t-\e,t) = \e^{-1}\int_0^1 \bigg( \bigg( I +
      \frac 13 \e^2 \lambda(t)^2 - \e^2\lambda(t)^2 y^2 \bigg)
       - \e^2 \sigma(0,t) y^2 \bigg) dy + \mathcal{O}(\e^2)
\\
   = \e^{-1}I-\frac \e 3 \sigma(0,t) + \mathcal{O}(\e^2) .  
\end{multline*}
Then, as $\rho(0,t)=\sigma(0,t)$, we obtain
\begin{equation*}
   {\bf S}(t-\e,t)^{-1}
       = \e \left(\mstrut{0.4cm}\right. I
               +\frac {\e^2} 3 \rho(0,t)
       + \mathcal{O}(\e^3) \left.\mstrut{0.4cm}\right) ,
\end{equation*}
so that
\[
   {\bf K}(r,t) = (t-r) \, I + \frac {(t-r)^3} 3 {\bf R}(t)
               + \mathcal{O}((t-r)^4)
\]
using that $\rho(0,t) = {\bf R}(t)$.
\end{proof}

\bigskip

\subsection{Reconstruction}

In this section we 
complete the proof of Theorem \ref{thm: main}. The main part is the proof of claim (i) which is given in the following proposition.
\medskip

\begin{proposition} \label{thm:1}
Functions ${\bf s}^j_k(0,t)$, $t > 0$, determine
uniquely functions ${\bf r}^j_k(r)$ for $r \in [0,t_0]$, and ${\bf s}^j_k(r,t)$ for $r,t \in [0,t_0]$ with $r<t$ where it is defined.
\end{proposition}

\bigskip

\begin{proof}
We are given the matrices ${\bf S}(0,t) = ({\bf
  s}^j_k(0,t))_{j,k=1}^{n-1}$, $t > 0$. Using
Lemma~\ref{lem:1}, it follows that the curvature matrix, ${\bf R}(r) =
({\bf r}^j_k(r))_{j,k=1}^{n-1}$, satisfies
\begin{equation} \label{eq:R3K}
   {\bf R}(r) = \frac 12 \p_t^3 {\bf K}(r,t)|_{t=r} ;
\end{equation}
similarly, ${\bf R}(r) = -\frac 12 \p_r^3 {\bf K}(r,t)|_{r=t}$.

Using
\[
   \p_r {\bf S}(r,t) = {\bf S}(r,t)^2 + {\bf R}(r)
\]
(cf.~(\ref{Matrix Riccati})) we find that
\begin{eqnarray*}
   \p_r {\bf K}(r,t) &=&
    -({\bf S}(r,t))^{-1} \p_r {\bf S}(r,t) ({\bf S}(r,t))^{-1}
\\
   &=& -({\bf S}(r,t))^{-1} ({\bf S}(r,t)^2
          + {\bf R}(r)) ({\bf S}(r,t))^{-1}
\\
   &=& -I - ({\bf S}(r,t))^{-1} {\bf R}(r)
            ({\bf S}(r,t))^{-1}
\\
   &=& -I - {\bf K}(r,t) {\bf R}(r) {\bf K}(r,t) .
\end{eqnarray*}
We let $\p_t$ act on the final equation above, and obtain
\begin{eqnarray*}
   \p_r ((\p_t {\bf K})(r,t))
          &=& \p_t (-I - {\bf K}(r,t) {\bf R}(r) {\bf K}(r,t))
\\
          &=& -((\p_t {\bf K})(r,t) {\bf R}(r) {\bf K}(r,t)
                + {\bf K}(r,t) {\bf R}(r) (\p_t {\bf K})(r,t)) .
\end{eqnarray*}
Computing the second and third $t$-derivatives in a similar manner,
and denoting $V = V(r,t) = (V^j(r,t))_{j=0}^3$, $V^j(r,t) = \p_t^j
{\bf K}(r,t)$ and ${\bf R} = {\bf R}(r)$, we obtain the equations
\begin{eqnarray}
   \p_r V^0 &=& -I - V^0 {\bf R} V^0 ,
\label{eq:Vsys1}
\\
   \p_r V^1 &=& -(V^1 {\bf R} V^0 + V^0 {\bf R} V^1 ),
\\
   \p_r V^2 &=& -(V^2 {\bf R} V^0 + V^0 {\bf R} V^2
                         + 2V^1 {\bf R} V^1) ,
\\
   \p_r V^3 &=& -(V^3 {\bf R} V^0 + V^0 {\bf R} V^3
                         + 3V^2 {\bf R} V^1 + 3V^1 {\bf R} V^2 ).
\label{eq:Vsys4}
\end{eqnarray}
Since ${\bf R}$ depends on $V^3$, this system is ``closed". We define
the operator $\T$ by
\begin{equation}
   (\T V)(r) = V^3(r,r)
\end{equation}
so that, with (\ref{eq:R3K}), ${\bf R} = {\bf R}(r) = \frac{1}{2} (\T
V)(r)$. Hence,
\begin{eqnarray}\label{EQ A}
   \p_r V^0 &=& -I - \frac 12 V^0(\T V)V^0 ,
\\
   \p_r V^1 &=&- \frac 12 (V^1(\T V)V^0 + V^0(\T V)V^1) , \nonumber
\\
   \p_r V^2 &=& -\frac 12 (V^2(\T V)V^0 + V^0(\T V)V^2\nonumber
                       + 2V^1(\T V)V^1) ,
\\
   \p_r V^3 &=& -\frac 12 (V^3(\T V)V^0 + V^0(\T V)V^3
                       + 3V^2(\T V)V^1 + 3V^1(\T V)V^2) .\nonumber
\end{eqnarray}
We write this as
\[
   \p_r V(r,t) = F(V(r,t),(\T V)(r)) ,
\]
where the map $F$ is a polynomial of its variables. We then introduce
\[
   \F :\ W(r,t) \mapsto F(W(r,t),(\T W)(r))
\]
so that the system (\ref{EQ A}) of nonlinear differential equations attains the
form
\begin{equation} \label{final diff. eq.}
   \p_r V(r,t) = (\F V)(r,t) .
\end{equation}
Assuming that we are given ${\bf S}(0,t)$ with $t > 0$, we know
the initial data
\begin{equation} \label{final init. data}
   V_0(t) = V(0,t) = (\p_t^j ({\bf S}(0,t))^{-1})_{j=0}^3 .
\end{equation}
We now address whether the initial value problem (\ref{final
  diff. eq.})-(\ref{final init. data}) has a unique solution.

Let us now assume that $t_1 > 0$ is such that the geodesic
$\gamma([0,t_1])$ has no focal points. This implies that the matrix ${\bf K}(r,t)$ is bounded on $(r,t) \in [0,t_1]^2$. Let $\mathcal{K}$ be a prior bound on the Riemannian
curvature, that is, $\| R \|_g \le \mathcal{K}$ on
$\gamma([0,t_1])$. The constant $\mathcal{K}$ depends thus on the
$C^3$-norm of the metric. 
Using \cite[Cor.~2.4]{Petersen-book}, which implies that
\[
{\bf S}(r,t) \ge
  \frac{\cos(\sqrt{\mathcal{K}} (t-r))}{\sin(\sqrt{\mathcal{K}}
    (t-r))} \, \mathcal{K} \, I
 \]
 boundedness is then guaranteed if $t - r \le \pi/(4
\sqrt{\mathcal{K}})$. Let $0 < t_2 < t_1$ and
\[
   Y_{t_2} = C([0,t_2]_r;C([0,t_2]_t;\R^{(n-1) \times (n-1)})^4)
\]
equipped with the norm
\begin{eqnarray*}
   \| V \|_{Y_{t_2}}  &:=& \sup_{r \in [0,t_2]}
                    \| V(r,.) \|_{C([0,t_2];\R^{(n-1) \times (n-1)})^4}\\
 & = &\sup_{(r,t) \in [0,t_2]^2} \max_{j \in \{0,\ldots,3\}}
                    \| V^j(r,t) \|_{\R^{(n-1) \times (n-1)}} .
\end{eqnarray*}
It is immediate that
\[
   |V^3(r,r)| \le 
   \sup_{(r,t) \in [0,t_2]^2}
                    \| V^3(r,t) \|_{\R^{(n-1) \times (n-1)}}
              \le \| V \|_{Y_{t_2}} .
\]
If $B_{t_2}(\mathcal{R}) \subset Y_{t_2}$ is the zero centered ball of
radius $\mathcal{R}\geq 1$ in $ Y_{t_2}$, because $\F$ contains no
differentiation, we find that
\[
   \F :\ \overline B_{t_2}(\mathcal{R}) \to Y_{t_2}
\]
is (locally) Lipschitz, with Lipschitz constant $L(\mathcal{R})$,
that does not depend on $t_2$.

We reformulate the differential equations in integral form, $H V = V$,
with
\[
   H :\ Y_{t_2} \to Y_{t_2} ,\quad
   (HW)(r,t) = V_0(t) + \int_0^r \F(W(r^\prime ,t)) dr^\prime  ,\quad
                   r,t \in [0,t_2] .
\]
Clearly, $H :\ \overline B_{t_2}(\mathcal{R}) \to Y_{t_2}$ is
(locally) Lipschitz, with Lipschitz constant $t_2 \,L(\mathcal{R})$.
For $H$ to be a contraction, we need that
\[
   t_2 \, L(\mathcal{R}) < 1 .
\]
To guarantee that $H(\overline B_{t_2}(\mathcal{R})) \subset \overline
B_{t_2}(\mathcal{R})$, we require that
\[
   \| V_0 \|_{C([0,t_1];\R^{(n-1) \times (n-1)})^4}
            + t_2 \, (1 + 4 \mathcal{R}^3) < \mathcal{R} .
\]
We choose
\[
   \mathcal{R} = 2 \, \| V_0 \|_{C([0,t_1];\R^{(n-1) \times (n-1)})^4} +1
\]
where the norm $\| V_0 \|_{C([0,t_1];\R^{(n-1) \times (n-1)})^4}$ can also
be bounded in terms of $\mathcal K$ using Gronwall's lemma. Indeed, using \cite[Cor.~2.4]{Petersen-book}, equations (\ref{eq:Vsys1})-(\ref{eq:Vsys4}), Lemma~\ref{lem:1}, and Gronwall' s lemma, we can obtain the following basic ``forward"  estimates:
\[
\| V^0(r,t) \|_{\R^{(n-1) \times (n-1)}} \lesssim \mathcal{K}^{-1},
\]
\[
\| V^1(r,t) \|_{\R^{(n-1) \times (n-1)}} \lesssim \mathrm{e}^{\tau_1/3},
\]
\[
\| V^2(r,t) \|_{\R^{(n-1) \times (n-1)}} \lesssim 2 \mathcal{K} \mathrm{e}^{4\tau_1},
\]
and
\[
\| V^3(r,t) \|_{\R^{(n-1) \times (n-1)}} \lesssim 8 \mathcal{K}^3 (1 + \tau_1 \mathrm{e}^{8 \tau_1})\mathrm{e}^{\tau_1}
\]
for $r, t \in [0,t_1]$ and $t - r \le \tau_1$, $\tau_1 = \tau_1(\mathcal{K}) = \pi/(4 \sqrt{\mathcal{K}})$; the
  maximum of these results is an estimate for $\| V_0
  \|_{C([0,t_1];\R^{(n-1) \times (n-1)})^4}$ in terms of $\mathcal{K}$.
Then we choose
\begin{equation} \label{EQ B}
   t_2 = \frac{1}{2} \min\left(
             \frac{\pi}{4 \sqrt{\mathcal{K}}},
             \frac{1}{L(\mathcal{R})},
             \frac{\mathcal{R}}{2 \, (1 + 4 \mathcal{R}^3)}\right),
\end{equation}
we see using the Banach fixed point theorem that $H$ has a unique fixed point in
$\overline B_{t_2}(\mathcal{R})$. Thus (\ref{final
  diff. eq.})-(\ref{final init. data}) has a unique solution $V \in
Y_{t_2}$.

{\mtext 
Let  ${\mathcal C}_r$ be the set of those $t\in [0,\infty)$ for
which $\gamma(r)$ and $\gamma(t)$ are conjugate points.
Recall that we assume that we are given
 the functions ${\bf s}^k_l(0,t)$ for all $t \in \R_+\setminus {\mathcal C}_0$ 
 and $l,k \in \{0,1,2,\dots,n\}$.
Thus we know ${\bf S}(0,t)$ and  ${\bf K}(0,t) ={\bf S}(0,t)^{-1}$, 
 for all $t \in \R_+\setminus {\mathcal C}_0$.
 Let us choose
$t_2 > 0$ so that (\ref{EQ B}) is valid. Then the equation (\ref{final
  diff. eq.}) has a solution $V(r,t)$, $(r,t)\in A_{t_2}:= \{(r,t)\in [0,t_2]^2, \ r\leq t\}$ with initial data
  $V(0,t)=(\p_t^j{\bf K}(0,t))_{j=0}^3$. Solving for $V(r,t)$ using equation
(\ref{final diff. eq.})  gives us also the curvature matrix ${\bf
  R}(r)$ for $r \in [0,t_2]$ by applying $\T$ (i.e. we compute $\mathbf{R}
= \frac12 \T V$). We will now switch notation replacing $t_2$ by $t_1$ intuitively indicating that for the first step we can reconstruct $\mathbf{R}$ from $0$ to $t_1$.

Next we do step by step reconstruction showing how to reconstruct $\mathbf{R}$ on the whole interval from $0$ to $t_0$ in steps.
First, we observe that
for $t\in \R_+$ outside of the discrete set ${\mathcal C}_r$ for
which $\gamma(r)$ and $\gamma(t)$ are conjugate points,
the matrix ${\bf S}(r,t)$ is {\mltext well defined.}

For given $t > t_1$, $t\not \in {\mathcal C}_0$, we will next reconstruct
${\bf S}(r,t)$, $r\in [0,t_1]$. As we know ${\bf R}(r)$ 
for $r\in [0,t_1]$ and the
 matrices ${\bf S}(0,t)$, we find on the interval $r \in
[0,t_1]$ the solutions ${\bf j}^k(r,t)$ of
the Cauchy
problems for  the Jacobi equations
\begin{eqnarray}
   \p_r^2 {\bf j}^k_l(r,t) + {\bf r}^k_p(r) {\bf j}^p_l(r,t)
                  &=& 0 ,\quad r \in [0,t_1] ,
\label{C-problems}
\\
   {\bf j}^k_l(r,t)|_{r=0} &=& \delta_l^k ,\quad
     \p_r {\bf j}^k_l(r,t)|_{r=0} = -{\bf s}_l^k(0,t).\nonumber
\end{eqnarray}
Now, when $t$ is given, for all $r\in [0, t_1]\setminus {\mathcal C}_t $
the vectors 
$\{{\bf  j}^k_l(r,t)\}_{l=1}^n$ are linearly independent; then, the equations $-\partial_r
{\bf j}^j_l(r,t) = {\bf s}_k^j(r,t) {\bf j}^k_l(r,t)$ (cf.~(\ref{Matrix
  shape})) determine ${\bf S}(r,t)$. Summarizing the above, 
  for all $t\in \R_+\setminus {\mathcal C}_0$ and
  $r\in [0,t_1]\setminus {\mathcal C}_t$ we can determine 
  ${\bf S}(r,t)$. As $t\mapsto {\bf S}(r,t)$ is continuous
  for $t\in  \R_+\setminus {\mathcal C}_r$, 
  we see that we can find ${\bf S}(r,t)$ for all
  $r\in [0,t_1]$ and $t>r$ such that 
$t \in \R_+\setminus {\mathcal C}_{r}$.
 In particular
we can determine ${\bf S}(t_1,t)$ for all $t>t_1$ such that 
$t \in \R_+\setminus {\mathcal C}_{t_1}$.
   This yields a
new dataset in the interior of $M$ at the point $\gamma(t_1)$. We now repeat the above argument with $0$ replaced by $t_1$ to recover $\mathbf{R}$ {\mltext and 
${\bf S}(r,t)$ on} another interval $[t_1,t_2]$.

As the metric is assumed to be $C^3$-smooth, the size of the
steps (i.e. $t_2-t_1$) are in a compact set 
uniformly bounded below by the right hand side of (\ref{EQ B}). Thus we can complete the reconstruction in a finite number
of steps and the proof is complete.}
This completes the proof of the claim (i) of Theorem \ref{thm: main}. 
\end{proof}

\section{Reconstruction of the metric in Fermi coordinates and 
the reconstruction the universal covering space}

{\mltext
We are now ready to prove claim (ii) of Theorem \ref{thm: main}.
\medskip

\begin{proof} 
Let us now recall some considerations from section \ref{sec:3}.
For any $x_0 \in U= \tilde M \setminus M$ and $\eta_0 \in T_{x_0} \tilde M$
and $t_0>0$ we considered the spherical surface $\Sigma_0\subset U$ with a center $y_{t_0}
=\gamma_{x_0,\eta_0}(t_0)$. On $\Sigma_0$
we considered the
 coordinates $\hat X:\hat V\to \R^{n-1}$ in a neighborhood 
$\hat V\subset \Sigma_0$ of $x_0$
and on the geodesics $\gamma_{\hat x}(r)$ we defined
the parallel frames $F_j(\hat x,r)$, $j=1,2,\dots,n$.
Note that as $\Sigma_0\subset U$ and we know the
metric tensor $g$ on $U$, we can determine the inner products
\begin{equation}\label{hat g matrix B}
   \hat g_{jk}(\hat x) = g(F_j(\hat x),F_k(\hat x)),
\end{equation}
and $g(F_j(\hat x,r),F_k(\hat x,r)) =\hat{g}_{jk}(\hat x)$ for all $r\geq 0$.

By the proof of claim (i) of Theorem \ref{thm: main},
for any $\hat x\in \hat X(\hat V)$  we can determine 
the coefficients $\mathbf{j}^k_l(r,t_0)=\mathbf{j}^k_l(\hat x,r,t_0)$ 
 given in (\ref{C-problems}). Then 
 \ba
J_j(\hat x,r;t_0)=\mathbf{j}_j^m(\hat x,r,t_0)F_m(\hat x,r)
 \ea
are the Jacobi fields along the geodesic $\gamma_{\hat x}(r)$
that satisfy
\ba
J_j(\hat x,t_0;t_0)=0,\quad J_j(\hat x,0;t_0)=F_m(\hat x,0).
\ea 

Let us now consider the set $W\subset \tilde M$ that
is a neighborhood of $\gamma_{x_0,\eta_0}((0,t_0)\setminus \mathcal C(x_0,\eta_0,t_0))$
and
the map $X_{\hat V}:\W\to W$ given in
 (\ref{eq: the map Xhat V}), a point $(\hat{x}_0,r_0)\in \W$
 and its small neighborhood $\V\subset \W$ so that $ X_{\hat V}|_\V:\V\to V=
 X_{\hat V}(\V)$ is a diffeomorphism. The inverse of this map defines local 
 coordinates  $x\mapsto (\hat{x},r)=(X_{\hat V}|_\V)^{-1}(x)$  in the set $V$.
We see that the Jacobi fields $J_j(\hat x,r)$ are in fact the coordinate vectors for the $(\hat{x},r)$ coordinates. Therefore the metric $g$ with respect to these coordinates can be recovered by
\[
g(\partial_{\hat{x}^j}, \partial_{\hat{x}^k}) |_{(\hat{x},r)} = \hat{g}_{ml}(\hat x)\mathbf{j}_j^m(\hat x,r;t_0) \mathbf{j}^l_k(\hat x,r;t_0) .
\]
Note that here  we are in fact varying $\hat{x}$, and performing the entire recovery of $\mathbf{r}$ and $\mathbf{s}$ along each geodesic in order to calculate the Jacobi fields along that geodesic
and can then compute the metric tensor $g$ in the set $W$ in the local $(\hat{x},r)$ coordinates.
Moreover, as we know the coefficients of the Jacobi fields  $J_j(\hat x,r)$ represented
in 
the parallel frame $F_m(\hat x,r)$ along $\gamma_{\hat x}(r)$,
{\mnewtext we can also  find the coefficients of the vectors $F_m(\hat x,r)$ in the basis given by the Jacobi fields $J_j(\hat x,r)$. Thus}
we change the local $(\hat{x},r)$ coordinates to 
the Fermi coordinates and  determine the metric tensor $g$ in the 
Fermi coordinates in some neighborhood $W^{fermi}_{x_0,\eta_0,t_0}\subset W$
of the set $\gamma_{x_0,\eta_0}((0,t_0)\setminus \mathcal C(x_0,\eta_0,t_0))$. 

As $W^{fermi}_{x_0,\eta_0,t_0}$ is a neighborhood of $\gamma_{x_0,\eta_0}((0,t_0)\setminus \mathcal C(x_0,\eta_0,t_0))$ we have not yet reconstructed $g$ in the whole
neighborhood of $\gamma_{x_0,\eta_0}$. To do this,
let $s_1>0$ be so small that $\tilde x_0=\gamma_{x_0,\eta_0}(-s_1)$ and
$\tilde \eta_0=\dot \gamma_{x_0,\eta_0}(-s_1)$  satisfy $\tilde x_0 \in U$
and repeat the above construction by replacing $x_0$ by $\tilde x_0$,
$\eta_0$ by $\tilde \eta_0$ and $t_0$ by
arbitrary $\tilde t_0>s_1$ and the spherical surface $\Sigma_0$ by the corresponding
surface $\tilde \Sigma_0$. Then, we can determine the metric
tensor in local coordinates $W^{fermi}_{\tilde x_0,\tilde \eta_0,\tilde t_0}$.
By varying $s_1$ and $\tilde t_0$ and using the fact that on a given geodesic the conjugate
points of a given point form a discrete set, we see that the whole geodesic 
$\gamma_{x_0,\eta_0}(\R_+)$ can be covered by neighborhoods of the form
$W^{fermi}_{\tilde x_0,\tilde \eta_0,\tilde t_0}$.
This completes the proof of claim (ii) of Theorem \ref{thm: main}.
\end{proof}
\medskip
}

\begin{figure}
\centering
\includegraphics[width=80mm]{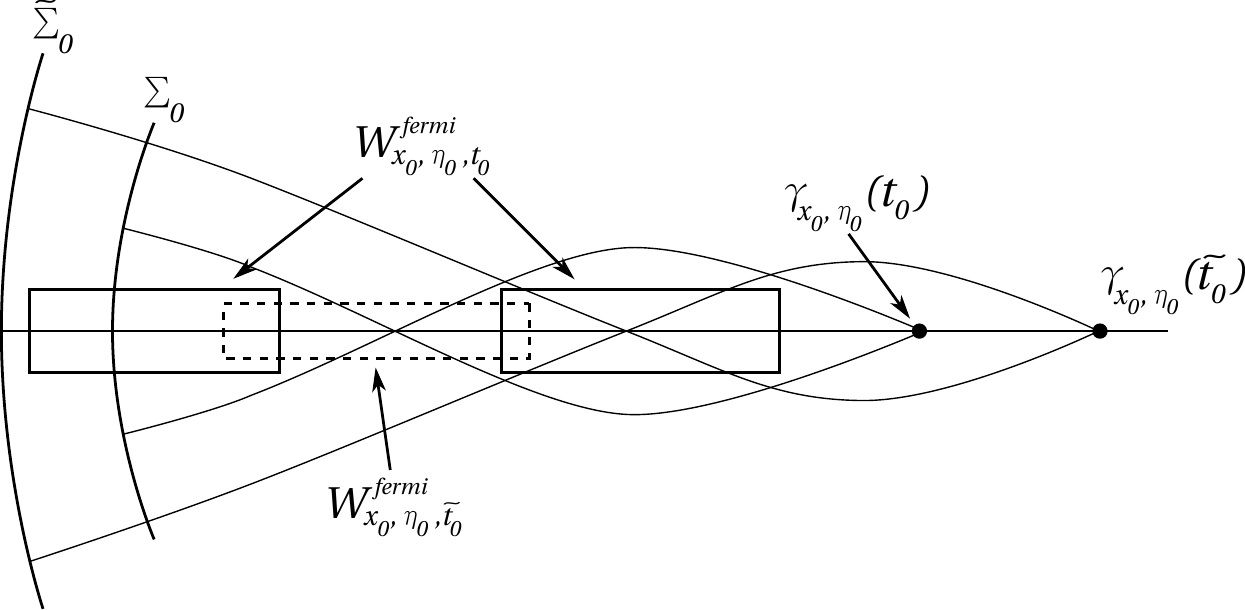}
\caption{Reconstruction procedure in the case of conjugate points.}
\label{fig:caustics}
\end{figure}


{\mltext
We finish this section by proving Theorem \ref{thm: main 2}. 
\medskip

\begin{proof} Let $\exp_x:T_x\tilde M\to \tilde M$ and  $\exp^\prime_{x^\prime}:T_{x^\prime}\tilde M'\to \tilde M^\prime$
be the exponential maps of $(\tilde M,g)$ and  $(\tilde M^\prime ,g^\prime )$, correspondingly. 

Let $p\in U$ and $p^\prime \in U^\prime $ be such
that $p^\prime =\Phi(p)$ and let
\ba
\ell=d\Phi|_p:T_p\tilde M\to T_{p^\prime }\tilde M^\prime 
\ea
be the differential of $\Phi$ at $p$. For $v=tv^0\in T_p\tilde M$, $\|v^0\|_g=1$, $t\geq 0$,
let $ \tau_v:T_p\tilde M\to T_q\tilde M$ denote
the parallel transport  along the geodesic $\gamma_{p,v^0}([0,t])$, where $q=\gamma_{p,v^0}(t)$
in $(\tilde M,g)$
and let $  \tau^\prime _{v^\prime }:T_{p^\prime }\tilde M^\prime \to T_{q^\prime }\tilde M^\prime $ denote the corresponding operation
on $(\tilde M^\prime ,g^\prime )$.  For $v,w\in T_p\tilde M$ let the curve
 $\mu_{v,w}:[0,2]\to \tilde M$ be the broken geodesic,
that is defined by
\ba
& &\mu_{v,w}(s)=\exp_p(sv),\quad\hbox{for $0\leq s\leq 1$,}\\
& &\mu_{v,w}(s)=\exp_q((s-1)\tau_v w),\quad\hbox{for $1\leq s\leq 2$,}
\ea
where $q=\exp_p(v)$. When $r=\mu_{v,w}(2)$ is the end point
of the broken geodesic, we denote by $\tau_{v,w}:
T_p\tilde M\to T_r\tilde M$ the parallel  transport of vectors
along the curve $\mu_{v,w}([0,2])$.
For all  $v\in T_p\tilde M$ let $\rho(p,v,r)$ be the function
in Theorem \ref{thm: main} for the geodesic $\gamma_{p,v}(r)$,
that is, we can determine  
the Riemannian metric
in the Fermi coordinates in
the tubular neighborhood $V_{x_0,\eta_0,r}$ that contains the ball
$B_{\tilde M}(\exp_p(v),\rho(p,v,r))$.
Let $\rho_0(v)>0$ be such that $\rho_0(v)<\rho(p,v,r)$ and the ball $B_{\tilde M}(\exp_p(v),\rho_0(v))$ is geodesically convex.
Let
$\rho^\prime _0 (v^\prime )>0$ be the corresponding function for $(\tilde M^\prime ,g^\prime )$ and $p^\prime $.
Finally, we define $f(v)=\min(\rho_0(v),\rho_0^\prime (\ell(v))).$

Let $v\in T_p\tilde M$. 
When $q=\exp_p(v)$ and $q^\prime =\exp^\prime _{p^\prime }(\ell(v))$, let $E_v:T_p\tilde M\to \tilde M$ 
be the map $E_v(\xi)=\exp_q( \tau_v\xi)$ and
$E^\prime _v:T_p\tilde M\to \tilde M$ be $E^\prime _{v^\prime }(\xi^\prime )=\exp^\prime _{q^\prime }( \tau^\prime _{v^\prime }\xi^\prime )$.
We see that \ba
E^\prime _{\ell(v)}\circ\ell\circ E_v^{-1}:B_{\tilde M}(q,f(v))\to 
B_{\tilde M^\prime }(q^\prime ,f(v))\ea is an isometry. In particular, if
 $v,w\in T_p\tilde M$ are such that $\|w\|_g<f(v)$, 
 and $r=\mu_{v,w}(2)\in \tilde M$ and $r^\prime =\mu^\prime _{\ell(v),\ell(w)}(2)\in \tilde M^\prime $
 are the end points of the broken geodesics,
 the above implies that the linear isometry 
 \ba
 \ell_{v,w}= \tau^\prime _{\ell(v),\ell(w)}\circ \ell \circ \tau_{v,w}^{-1}:
 T_r\tilde M\to  T_{r^\prime }\tilde M^\prime
  \ea
 preserves the sectional curvature, i.e.,
 $\hbox{Sec}_{r}(\xi,\eta)=\hbox{Sec}_{r^\prime }(\ell_{v,w}(\xi),\ell_{v,w}(\eta))$.
 Thus from the proof of the Ambrose theorem given  in \cite{ONeill-Proc} (for the original
 reference, see \cite{Ab}) it
 follows that $(\tilde M,g)$ and $(\tilde M^\prime ,g^\prime )$ have isometric covering spaces.
  \end{proof}
 }

\medskip\medskip


\def\cprime{$^\prime $} \newcommand{\SortNoop}[1]{}

\section*{Appendix A: Proof of Proposition \ref{prop: determination in U}}
\label{subsec: Equivalence}


We begin with the reconstruction of the Riemannian metric $g|_U$ if
$\mathcal S_U$ is given. If $x \in U$ then in a sufficiently small
neighborhood $U^\prime  \subset U$ of $x$ all points $z$ can be connected to
$x$ with a geodesic of a given length (travel time) contained in
$U$. As a consequence, the distance between $x$ and $z$ can be found
as 
\begin{multline}
   d(x,z) = \inf \left\{\mstrut{0.5cm}\right.
   \sum_{j=1}^N 2t_j ;\ 
                     (t_j,\Sigma_j) \in {\mathcal S}_U ,\
                     x_j,x_{j+1} \in \Sigma_j
\\[-0.3cm]
   \hbox{for $j=1,\dots,N$ such that $x_1 = x$, $x_{N+1} = z$}
                 \left.\mstrut{0.5cm}\right\} .
\end{multline}
Indeed, we observe that the infimum is obtained when $\Sigma_j$ are
the boundaries of sufficiently small balls (which are always smooth),
$B_{\tilde M}(y_j,t_j)$, where $y_j$ are points on the shortest geodesic
connecting $x$ to $z$. Thus we can determine the distance function
$(y,y^\prime ) \mapsto d(y,y^\prime )$ between two arbitrary points in $y, y^\prime  \in
U^\prime $.

Now, if $r > 0$ is small enough and $z_1,\dots,z_n \in U^\prime $ are
disjoint points so that $d(x,z_j) = r$, then the function $y \mapsto
(d(y,z_j))_{j=1}^n \in \R^n$ defines local coordinates near the point
$x \in U^\prime $. So, in $U^\prime $, we can find the differentiable structure inherited from
the manifold $\tilde M$. Using the distances beween points $y, y^\prime  \in
U^\prime $, we can determine the Riemannian metric in these coordinates in
$U^\prime $. But then we can find the Riemannian metric $g|_U$ if ${\mathcal
  S}_U$ is given.

{\mnewtext
For $(t,\Sigma) \in {\mathcal S}_U$ and
$x_0 \in \Sigma$, let $N(x_0,\Sigma,t)$ be the set consisting
of the two unit normal vectors of $\Sigma$ at $x_0$.  Let $ N_1(x_0,\Sigma,t)$
be the set of those $\nu_0\in N(x_0,\Sigma,t)$ for which the point $x_0$ has a neighborhood $U^\prime 
\subset U$ such that $ \Sigma \cap U^\prime$ has the representation
\begin{equation} \label{Sigma representation}
   \Sigma \cap U^\prime  = \{\gamma_{y,\eta}(t) ;\
       \eta \in \Omega \} ,
\end{equation}
where $y = \gamma_{x_0,\nu_0}(-t)$ and $\Omega \subset \Omega_y \tilde M$
is a neighborhood of $\eta_0 = -\dot\gamma_{x_0,\nu_0}(-t)$. 
Note that it is possible for $N_1(x_0,\Sigma,t)$ to contain both
normal vectors in $N(x_0,\Sigma,t)$. An
example of a case in which this occurs is  when
 $\tilde M$ is $S^2$ and $\Sigma$ is a subset of the Equator.

\medskip\medskip

\begin{lemma} \label{lem: normal vector}  If $U$ and
${\mathcal S}_U$ are given, we can determine $N_1(x_0,\Sigma,t)$ for any $(t,\Sigma) \in {\mathcal S}_U$ and
$x_0 \in \Sigma$.
\end{lemma}
\medskip


\begin{proof} 
For given $(t,\Sigma)\in \mathcal{S}_U$ and $x_0\in U$
let  $\zeta_0\in N(x_0,\Sigma,t)$ be one of the two unit normal vectors to $\Sigma$ at $x_0$, and let
$\zeta(x)$ be a smooth normal vector field on $\Sigma \cap
U^\prime $ such that $\zeta(x_0) = \zeta_0$. We introduce the notation
\[
   \Sigma^{\pm}_{U^\prime }(s)
        = \{ \gamma_{x,\pm\zeta(x)}(s) ;\ x \in \Sigma \cap U^\prime  \} ;
\]
$\Sigma^{\pm}_{U'}(s)$ will be smooth for $s \in (-\e,\e)$, $\e > 0$
sufficiently small.

Assume next that  $\zeta_0 \in N_1(x_0,\Sigma,t)$. 
Then representation (\ref{Sigma representation}) is
valid with $y=\gamma_{x_0,\zeta(x_0)}(-t)$, and for $p(x,s) =
\gamma_{x,\zeta(x)}(s)$ and $\eta(x,s) = \dot\gamma_{x,\zeta(x)}(s)$,
we have 
$\gamma_{p(x,s),\eta(x,s)}(-t-s) = \gamma_{x,\zeta(x)}(s-t-s)= y$.  Hence,
there is a neighborhood $U^{\prime\prime} \subset U$ of $x_0$ such that for all $\e>0$ small enough
\beq\label{eq: plus condition}
   (t+s,\Sigma^+_{U^{\prime\prime} }(s)) \in {\mathcal S}_{U }
               \hbox{ for all }s \in (-\e,+\e) .
\eeq
Let us consider the following condition
\ba
\hbox{(C) There exists 
 $y^\prime $ such that $\gamma_{x,\zeta(x)}(+t) = y^\prime $ for all
$x \in \Sigma$ close to $x_0$.}
\ea
If condition (C) is valid, then both $\zeta_0$ and $-\zeta_0$
are in $N(x_0,\Sigma,t)$. If  condition (C) is not valid, then 
$N(x_0,\Sigma,t)$ contains the vector $\zeta_0$ but not $-\zeta_0$.

In the case when condition (C) is valid, we see that (\ref{eq: plus condition}) holds
as well as the analogous identity with the minus sign, that is, we have
\beq\label{eq: minus condition}
   (t+s,\Sigma^-_{U^{\prime\prime} }(s)) \in {\mathcal S}_{U^{\prime\prime} }
               \hbox{ for all }s \in (-\e,+\e) .
\eeq
Next, consider the case when the condition (C) is not valid.
Our aim is show that then (\ref{eq: minus condition}) can not hold.
For this end, let us assume that the condition (C) is not valid
but we have  (\ref{eq: minus condition}). Then we see that
 for all $s\in (-\e,\e)$ one of the sets \ba
& &   A_{+} = \{ \gamma_{x,\zeta(x)}(-s + (t+s)) ;\
                  x \in  \Sigma \cap U^{\prime\prime}  \},\\
& &   A_{-} = \{ \gamma_{x,\zeta(x)}(-s - (t+s)) ;\
                  x \in \Sigma \cap U^{\prime\prime}  \} ,
\ea
would consist of a single point.  Now, if $A_+$ consisted of a single
point then this point would satisfy the condition required for the
point $y^\prime $ in condition (C).  As we assumed that 
condition (C) is not valid, we
conclude that $A_+$ cannot be a single point. If $A_-$ consisted of a single point for all $s \in
(-\e,\e)$, then for all $x_1,x_2 \in \Sigma\cap U^{\prime\prime} $ we would
have $\gamma_{x_1,\zeta(x_1)}(-t-2s)=\gamma_{x_2,\zeta(x_2)}(-t-2s)$ for all $s\in (-\e,\e)$ and hence for all $s\in \R$. With $s=-t/2$
we would see that $x_1=x_2$ for all $x_1,x_2 \in \Sigma\cap U^{\prime\prime}$ but that is not possible. Hence  equation (\ref{eq: minus condition}) can not
be true when the condition (C) is not valid

Summarizing the above, we can find the set $N_1(x_0,\Sigma,t)$ using the fact that it contains the vector $\pm\zeta_0$ if 
and only if there are $U^{\prime\prime} $ and $\e>0$ such that
$(t+s,\Sigma^\pm_{U^{\prime\prime}}(s)) \in {\mathcal S}_{U}$ holds for all $s \in (-\e,\e)$. 
\end{proof}
}

\medskip
Lemma \ref{lem: normal vector} and the considerations above it prove Proposition \ref{prop: determination in U}.

\section*{Bibliography}

\bibliographystyle{model1-num-names}

\end{document}